\documentclass[12pt,a4paper]{article}
\usepackage{latexsym}
\usepackage{amsfonts, dsfont, eufrak, yfonts}
\usepackage{amssymb,fontenc,mathrsfs,bbm}
\usepackage[notref,notcite]{showkeys}
\newcommand{\grad}{{{\rm{grad}}}}

\newcommand{\Hess}{{{\rm Hess}}}
\newcommand{\noin}{\noindent}

\newtheorem{theorem}{\bf Theorem}[section]

\newtheorem{corollary}[theorem]{\bf Corollary}

\newtheorem{remark}[theorem]{\bf Remark}

\begin{document}

\title{On cylindrically bounded $H$-Hypersurfaces of $\mathbb{H}^{n}\times \mathbb{R}$}
\author{G. Pacelli Bessa\thanks{Research partially supported by CNPq-Brasil} \and M. Silvana Costa}
\date{\today}
\maketitle
\begin{abstract}
\noin We show that $H$-hypersurfaces of $\mathbb{H}^{n}\times
\mathbb{R}$ contained in a vertical cylinder and with Ricci
curvature with strong quadratic decay have mean curvature $\vert
H\vert > (n-1)/n$.\\


\end{abstract}

\section{Introduction}Barbosa-Kenmotsu-Oshikiri in
 \cite{barbosa-kenmotsu-Oshikiri} and Salavessa in  \cite{salavessa1} independently
  proved that if an entire graph of a smooth function
$ f:\mathbb{H}^{n}\to \mathbb{R}$ has constant mean curvature $H$
then $\vert H \vert\leq (n-1)/n$. On the other hand,
Bessa-Montenegro in \cite{bessa-montenegro} showed that the mean
curvature of any compact $H$-hypersurface  immersed in
$\mathbb{H}^{n}\times \mathbb{R}$  satisfies $\vert H \vert >
(n-1)/n$.  It should be remarked that  this result for embedded
$H$-surfaces ($n=2$) was proved by Nelli-Rosenberg in \cite{nelli}
and it was implicit in Hsiang-Hsiang's paper \cite{hsiang-hsiang}.
The case of embedded $H$-hypersurfaces ($n\geq 2$) follows from
Salavessa's work \cite{salavessa1}, \cite{salavessa2}.  The purpose
of this paper is to extend Bessa-Montenegro's result to
cylindrically bounded complete hypersurfaces $M$ of
$\mathbb{H}^{n}\times \mathbb{R}$ with Ricci curvature  with strong
quadratic decay. We say that a complete Riemannian manifold $M$ has
Ricci curvature $Ric_{M}$ with strong quadratic decay if
$$Ric_{M}(x) \geq -c^{\,2}\left[1+ \rho^{2}_{M}(x)\log^{2}
(\rho_{M}(x)+2)\right],$$ where $\rho_{M}$ is the distance function
on $M$ to a fixed point $x_{0}$ and $c=c(x_{0})>0$ is a constant
depending on $x_{0}$. An immersed submanifold $M\subset N\times
\mathbb{R}$ is said to be cylindrically bounded if $p(M)$ is a
bounded subset of $N$, where $p:N\times \mathbb{R}\to N$ is the
projection on the first factor $p(x,t)=x$. Our main result is the
following theorem.
\begin{theorem}\label{thm1}Let $M$ be a complete hypersurface immersed in
 $\mathbb{H}^{n}\times
\mathbb{R}$ with Ricci curvature  with strong quadratic decay. If
$M$ is cylindrically bounded then $\sup_{M}\vert H \vert \geq
(n-1)/n$.
\end{theorem} Calabi, \cite{calabi} in the sixties asked whether
there were complete bounded minimal surfaces in $\mathbb{R}^{3}$. It
is well known that this question was completely answered. See
\cite{omori}, \cite{jorge-xavier}  for non-existence of complete
bounded minimal surfaces with bounded sectional curvature, see
\cite{chen-xin} for non-existence of complete bounded minimal
surfaces   with sectional curvature with strong quadratic decay and
see \cite{Nadirashvili} for the first example of a complete bounded
minimal surface in $\mathbb{R}^{3}$.  Naturally, one can ask whether
there are complete bounded  minimal surfaces in
$\mathbb{H}^{2}\times \mathbb{R}$. The answer is not if the minimal
surface has sectional curvature with  strong quadratic decay. This
is implied by  Theorem (\ref{thm1}). In fact, it  can be restated as
the following corollary.

\begin{corollary}Let $M$ be a complete $H$-hypersurface immersed in
 $\mathbb{H}^{n}\times
\mathbb{R}$ with Ricci curvature  with strong quadratic decay. If
$\vert H \vert < (n-1)/n$ then $M$ is not cylindrically bounded.
\end{corollary}
It is well known that the coordinate functions of any minimal
surface of $\mathbb{R}^{3}$ are harmonic, thus if one of the
coordinates  is bounded then the minimal surface is non-parabolic.
In our second result,  we show that cylindrically bounded
$H$-hypersurfaces of $\mathbb{H}^{n}\times \mathbb{R}$ with $\vert
H\vert < (n-1)/n$ are non-parabolic.
\begin{theorem}\label{thm2}Let
$M$ be a complete immersed hypersurface in  $\mathbb{H}^{n}\times
\mathbb{R}$ with bounded mean curvature $\sup_{M}\vert H \vert \leq
(n-1)/n$ contained in a vertical cylinder. Then $M$ is
non-parabolic.
\end{theorem}Any hypersurface of $\mathbb{S}^{n}\times \mathbb{R}$ is
cylindrically bounded, nevertheless, the Theorems (\ref{thm1}) and
(\ref{thm2}) have appropriate versions for $\mathbb{S}^{n}\times
\mathbb{R}$.

\begin{theorem}\label{thm3}Let $M$ be a complete hypersurface immersed in
 $\mathbb{S}^{n}\times
\mathbb{R}$ with Ricci curvature  with strong quadratic decay. If
$p_{1}(M)\subset B_{\mathbb{S}^{n}}(r)$, $r< \pi/2$  then
$\sup_{M}\vert H \vert \geq (n-1)\cot(r)/n$.
\end{theorem}

\begin{theorem}\label{thm4}Let $M$ be a complete hypersurface immersed in
 $\mathbb{S}^{n}\times
\mathbb{R}$ with bounded mean curvature $\sup_{M}\vert H\vert \leq
(n-1)\cot(r)/n$, $r\leq \pi/2$. Then $M$ is non-parabolic.
\end{theorem}

\section{Preliminaries}
Let $\varphi : M \hookrightarrow N$ be an isometric immersion, where
$M$ and $N$ are complete Riemannian
 $m$ and $n$ manifolds respectively. Consider a smooth function $g:N \rightarrow \mathbb{R}$ and the composition $f=g\,\circ\,
 \varphi :M \rightarrow \mathbb{R}$. Identifying $X$ with $d\varphi (X)$ we have at $q\in M$ and  for every
 $X\in T_{q}M$ that
\begin{eqnarray} \langle \grad\,f\, , \, X\rangle \,=\,d f(X)=d g(X)\,=\, \langle \grad\,g\, ,\, X\rangle,\nonumber
\end{eqnarray}
therefore
\begin{equation} \grad\, g\,=\, \grad\, f\, +\,(\grad\,g)^{\perp},\label{eqBF1}
\end{equation}
where $(\grad\,g)^{\perp}$ is perpendicular to $T_{q}M$. Let $\nabla
$ and $\overline{\nabla}$ be the Riemannian connections on $M$ and
$N$ respectively,  $\alpha (q) (X,Y) $ and $\Hess\,f(q)\,(X,Y)$  be
respectively the second fundamental form of the immersion $\varphi $
and
 the Hessian of $f$ at $q\in M$,  $X,Y \in T_{p}M$. Using the Gauss equation we have that
 \begin{equation}\Hess\,f (q)  \,(X,Y)= \Hess\,g (\varphi (q))\,(X,Y) +
 \langle \grad\,g\,,\,\alpha (X,Y)\rangle_{\varphi (q)}.
\label{eqBF2}
\end{equation}
 Taking the trace in (\ref{eqBF2}), with respect to an orthonormal basis $\{ e_{1},\ldots e_{m}\}$
 for $T_{q}M$, we have  that
\begin{eqnarray}
\Delta \,f (q) & = & \sum_{i=1}^{m}\Hess\,g (\varphi
(q))\,(e_{i},e_{i}) + \langle \grad\,g\,,\,
                             \sum_{i=1}^{m}\alpha (e_{i},e_{i})\rangle.\label{eqBF3}
\end{eqnarray}
These formulas (\ref{eqBF2}) and (\ref{eqBF3}) are  well known in
the literature,
 see \cite{bessa-montenegro1}, \cite{bessa-jorge-lima-montenegro},
  \cite{choe-gulliver}, \cite{cheung-leung}, \cite{jorge-koutrofiotis}.
  Recall the Hessian
  Comparison Theorem.

\begin{theorem} Let $M$ be a complete Riemannian manifold  and
$x_{0},x_{1} \in M $.   Let $\rho (x)$
 be the  distance function $dist_{M}(x_{0}, x) $ to $x_{0}$ and let
 $\gamma$ be a minimizing geodesic joining $x_{0}$ and $x_{1}$. Let $K_{\gamma}$ be the radial
  sectional curvatures of $M$ along $\gamma$
   and  let $\mu(\rho )$
  be this function defined below.

\begin{equation} \mu(\rho )=\left\{ \begin{array}{lcll}
 & k \cdot\coth (k \cdot\rho (x)), & if  & \sup K_{\gamma}=-k^{2} \\
 &                       &     &   \\
 &\displaystyle \frac{1}{\rho (x)},   & if  &  \sup K_{\gamma}=0  \\
 &                       &     & \\
 & k\cdot\cot (k\cdot \rho (x)),  & if  &  \sup K_{\gamma} =k^{2}\; and \; \rho < \pi/2k.
\end{array}\right.\label{eqBF5}
\end{equation}

\noin Then the Hessian of  $\rho$ satisfies

\begin{equation}\begin{array}{llllllll}

Hess\,\rho(x)(X,X)&\geq & \mu(\rho(x))\cdot\Vert
X\Vert^{2},&&Hess\,\rho(x) (\gamma ',\gamma ') & =   & 0
\end{array}\label{eqBF6}
\end{equation}
Where $X$ is any vector in $T_{x}M$ perpendicular to $\gamma'(\rho
(x))$. \label{thmHess}

\end{theorem}
The second main ingredient to the proof of our results is the
Omori-Yau maximum principle \cite{omori}, \cite{yau}, \cite{dias} in
the generalized version proved by  Chen-Xin in \cite{chen-xin}.

\begin{theorem}[Omori-Yau Maximum Principle, \cite{chen-xin}] Let
$M$ be a complete Riemannian manifold with Ricci curvature with
strong quadratic decay. Let $u$ be a $C^{2}$ function bounded above
on $M$. Then for any sequence $\epsilon_{k}\to 0$ of positive
numbers there exists a sequence of points $x_{k}\in M$ such that
\begin{itemize}
\item[i.] $\lim_{k\to \infty} u(x_{k})=\sup_{M}u$
\item[ii.] $ \vert \grad u\vert (x_{k})< \epsilon_{k}$
\item[iii.] $ \triangle u (x_{k})< \epsilon_{k}$.
\end{itemize}
\end{theorem}
\begin{remark}It is clear in the proof of  the Omori-Yau maximum principle that $u$ is allowed to be only $C^{0}$
 in a measure zero subset of $M$. So that $u$ can be the distance function on $M$ to a fixed point.
  See page 360 of \cite{chen-xin}.
\end{remark}

\section{Proof of the Results}\label{sectionProof}Our results (Thms. (\ref{thm1}, \ref{thm2}, \ref{thm3}, \ref{thm4})
 are particular case of a
more general theorem that we present here. First, let $\varphi : M
\hookrightarrow N\times \mathbb{R}$ be an isometric immersion of a
complete Riemannian $m$-manifold $M$  into a complete Riemannian
$n$-manifold $N$ with a pole and radial sectional curvature
$K_{N}\leq -\kappa^{2}<0$. Let $\rho_{N}:N\to \mathbb{R}$ be the
distance function  to the pole $x_{0}\in N$. Set $g: N\times
\mathbb{R}\to \mathbb{R}$ and $f=g\circ \varphi:M\to \mathbb{R}$.
The Laplacian of $f$ is given by
$$
\Delta \,f (q)    =   \sum_{i=1}^{m}\Hess\,g (\varphi
(q))\,(e_{i},e_{i}) + \langle \grad\,g\,,\,
                            \stackrel{\rightarrow}{H}\rangle.
                             $$Here
$\stackrel{\rightarrow}{H}$ is the mean curvature vector with norm
$\Vert \stackrel{\rightarrow}{H}\Vert =n \vert H \vert$.
 We choose   an orthonormal basis  $\{
e_{1},\ldots e_{n}\}$
 for $T_{q}M$ in the following way.
  Start with an orthonormal basis (from polar coordinates) for $T_{p_{1}(q)}N$,
  $\{\grad \rho_{N},\partial/\partial \theta_{2},\ldots,
  \partial/\partial \theta_{n}\}$. We can choose be an orthonormal  basis for $T_{q}M$ as
  follows
  $e_{1}=\langle e_{1}, \partial /\partial t\rangle \partial /\partial t
 + \langle e_{1}, \grad \rho_{N}\rangle\grad \rho_{N}$ and $e_{j}=\partial/\partial \theta_{j}$, $j=2,...,n$
 (up to an re-ordination), where $\partial/\partial t$ is tangent to the $\mathbb{R}$-factor.
  By the Hessian comparison theorem we
have $q\in M$ that \begin{equation}\label{Algazarra}\Delta \,f (q)
\geq (n-1)\kappa\cdot\coth (\kappa \cdot\rho_{N})) - n\vert H \vert
.\end{equation}Where the right hand side of the inequality
(\ref{Algazarra}) is computed at $\varphi (q)$.
\vspace{2mm}

Suppose that  $\varphi (M)$  has  mean curvature vector $ \Vert
\stackrel{\rightarrow}{H}\Vert=n \vert H \vert < (n-1)\kappa$. This
implies that the function $f$ is subharmonic. If $\varphi (M)$ is
cylindrically bounded then $f$ is a bounded subharmonic function and
$M$ is non-parabolic. Suppose in addition the $M$ has Ricci
curvature with strong quadratic decay.  By the Omori-Yau Maximum
principle there exist sequences $\epsilon_{k}\to 0$ and $x_{k}\in M$
with $\triangle f (x_{k})< \epsilon_{k}$. The inequality
(\ref{Algazarra}) at $x_{k}$  becomes
\begin{equation}\label{voraz}\epsilon_{k}>\triangle f (x_{k})\geq
 (n-1)\kappa-
n\sup_{M}\vert H \vert >0
\end{equation}Letting $\epsilon_{k}\to 0$ we get a contradiction. This
proves the following more general result from which Theorems
(\ref{thm1} and \ref{thm2}) are corollaries.

\begin{theorem}\label{thm5} Let $N$ be a complete Riemannian $n$-manifold with a pole and radial
sectional curvature bounded above, $K_{N}\leq -\kappa^{2}<0$. Let
$\varphi :M\hookrightarrow N\times \mathbb{R}$ be a complete
 immersed submanifold. Then
 \begin{itemize}
 \item[i.]If $\varphi (M)$ is  cylindrically bounded and  has bounded mean curvature vector $
\Vert \stackrel{\rightarrow}{H}\Vert < (n-1)\kappa$ then $M$ is
non-parabolic.
 \item[ii.]If $\varphi (M)$ is  cylindrically bounded and $M$ has Ricci curvature strong quadratic
decay then then $\sup_{M}\vert H \vert \geq (n-1)\kappa/n$.
 \end{itemize}
\end{theorem}

Now suppose that $\varphi :M\hookrightarrow \mathbb{S}^{n}\times
\mathbb{R}$  is a complete immersed submanifold. The inequality
(\ref{Algazarra}) becomes in this setting the following inequality.
For any $q\in M$
\begin{equation}\label{gigantesco}\Delta \,f (q)
\geq \left[(n-1)\cdot\cot (\rho_{\mathbb{S}^{n}}(\cdot )) - n\vert H
\vert\right](\varphi (q)) .\end{equation} If $p(\varphi (M))\subset
B_{\mathbb{S}^{n}}(r)$, $r< \pi/2$ then $\cot (r)>0$. Thus if
$\sup_{M}\vert H \vert < (n-1)\cot (r)/n$ then $f$ is a bounded
subharmonic function. And if in addition,  $M$ has Ricci curvature
with strong quadratic decay then we have by Omori-Yau maximum
principle that $$\epsilon_{k}>\triangle f (x_{k})\geq (n-1)\cdot\cot
(r) - n\sup_{M}\vert H \vert = b^{2}>0 $$ Letting $\epsilon_{k}\to
0$ we get a contradiction. This proves the following theorem that
implies Theorems (\ref{thm3} and \ref{thm4}) as  consequences.
\begin{theorem}\label{thm6}Let $\varphi :M\hookrightarrow \mathbb{S}^{n}\times
\mathbb{R}$  be a complete immersed submanifold.
\begin{itemize}
 \item[i.]If $p(\varphi (M))\subset
B_{\mathbb{S}^{n}}(r)$, $r< \pi/2$ and  $\sup_{M}\vert H \vert <
(n-1)\cot (r)/n$ then $M$ is non-parabolic.
 \item[ii.]If $p(\varphi (M))\subset
B_{\mathbb{S}^{n}}(r)$, $r< \pi/2$ and $M$ has Ricci curvature
strong quadratic decay then then $\sup_{M}\vert H \vert \geq
(n-1)\cot (r)/n$.
 \end{itemize}
\end{theorem}

\vspace{1cm}

\noin {\it Address of the authors:}\\
\noin Departamento de Matematica\\
\noin Campus do Pici, Bloco 914\\
\noin Universidade Federal do Cear\'{a}-UFC\\
\noin 60455-760 Fortaleza-Cear\'{a}\\
\noin Brazil\\

\noin  bessa@math.ufc.br \& mscosta@mat.ufc.br
\end{document}